\documentclass[fleqn]{mat01}
\usepackage{times,mathtimy,amssymb,epsfig}
\begin{document}

\setcounter{page}{107}
\firstpage{107}

\newtheorem{theore}{Theorem}
\renewcommand\thetheore{\arabic{section}.\arabic{theore}}
\newtheorem{theor}[theore]{\bf Theorem}
\def\notation{\trivlist\item[\hskip\labelsep{{\it Notation.}}]}
\def\teo{\trivlist\item[\hskip\labelsep{\bf Theorem}]}
\def\potc{\trivlist\item[\hskip\labelsep{{\it Proof of the Claim.}}]}
\newtheorem{theorr}{Theorem}
\renewcommand\thetheorr{\arabic{theorr}}
\newtheorem{therr}[theorr]{\bf Theorem}
\newtheorem{lem}[theore]{Lemma}
\newtheorem{propo}[theore]{\rm PROPOSITION}
\newtheorem{rem}[theore]{Remark}
\newtheorem{caa}{Case}
\renewcommand\thecaa{(\alph{caa})}
\newtheorem{case}[caa]{Case}
\newtheorem{coro}[theore]{\rm COROLLARY}
\newcommand{\A}{{\mathbb A}}
\newcommand{\C}{{\mathbb C}}
\newcommand{\V}{{\mathbb V}}
\newcommand{\Z}{{\mathbb Z}}
\def\o{{\cal O}}
\def\P{{\Bbb P}}
\def\L{{\Bbb L}}
\def\G{{\cal G}}
\def\E{{\cal E}}
\def\D{{\cal D}}

\font\zzzz=tibi at 10.4pt
\def\U{\mbox{\zzzz{U}}}
\def\n{\mbox{\zzzz{n}}}
\def\d{\mbox{\zzzz{d}}}
\def\U{\mbox{\zzzz{U}}}

\font\zzzzz=tibi at 7.3pt
\def\La{\mbox{\zzzzz{L}}}

\font\zzzzzz=tib at 7.3pt
\def\co{\mbox{\zzzzzz{,}}}

\title{Picard groups of the moduli spaces of semistable sheaves I}

\markboth{Usha N~Bhosle}{Pic of moduli of semistable sheaves}

\author{USHA N~BHOSLE}

\address{Tata Institute of Fundamental Research, Homi Bhabha Road,
Mumbai~400~005, India\\
\noindent E-mail: usha@math.tifr.res.in}

\volume{114}

\mon{May}

\parts{2}

\Date{MS received 9 January 2003; revised 12 March 2004}

\begin{abstract}
We compute the Picard group of the moduli space $U'$ of semistable
vector bundles of rank $n$ and degree $d$ on an irreducible nodal curve
$Y$ and show that $U'$ is locally factorial. We determine the canonical
line bundles of $U'$ and $U'_L$, the subvariety consisting of vector
bundles with a fixed determinant. For rank $2$, we compute the Picard
group of other strata in the compactification of $U'$.
\end{abstract}


\maketitle

\section{Introduction}

In our previous paper \cite{3} we proved that the Picard group of the
moduli space $U'_L(n,d)$ of semistable vector bundles of rank $n$ with
fixed determinant $L$ ($L$ being a line bundle of degree $d$) on an
irreducible projective nodal curve $Y$ of geometric genus $g \geq 2$ is
isomorphic to $\Z$ (except possibly in the case $g=2, n=2, d$ even). We
used this to show that $U'_L(n,d)$ is locally factorial. Interestingly,
the results for irreducible nodal curves are very similar to those for
smooth curves. However, the proofs are different and much more
difficult. Unlike in the smooth case, the moduli space of vector bundles
on a nodal curve is not projective. Moreover its complement in the
compactification $U$ (moduli of torsion-free sheaves) has codimension
$1$. The computation of Picard group needs codimension of the
non-semistable and non-stable strata (see \cite{6,11} for smooth case).
Since HN-filtrations of vector bundles contain non-locally free sheaves
and tensor products of stable bundles are not semistable (on $Y$), in
general it is impossible to determine this codimension directly on $Y$.
We did it by using parabolic bundles on the normalization $X$ of $Y$ and
hence had to assume $g \ge 2$ and exclude the case $g=n=d=2$.

In this paper, we do a detailed analysis for rank $2$ and extend these
results to nodal curves of arithmetic genus $g_Y \ge 0$ (rank $2$).
Combining this with results of \cite{3}, we have the following theorem.

\begin{therr}[\!]
Let $Y$ be an irreducible reduced curve with only ordinary nodes as
singularities. Assume that for $n \ge 3${\rm ,} the geometric genus $g
\ge 2$. Then
\begin{enumerate}
\renewcommand\labelenumi{{\rm (\arabic{enumi})}}
\item {\rm Pic} $U'_L(n,d) \approx$ {\rm Pic} $U_L^{\prime s}(n,d)
\approx \Z${\rm ,}\vspace{.15pc}

\item $U'_L$  is locally  factorial.
\end{enumerate}
\end{therr}

We also show that the dualising sheaf $\omega_L$ of $U'_L(n,d)$ is
isomorphic to the line bundle $-2 \delta \L$, where $\delta=gcd (n,d)$
and $\L$ is the ample generator of Pic $U'_L(n,d)$\break (Theorem 4).

We then compute the Picard group of the moduli space $U'(n,d)$ (resp.
$U^{\prime s} (n,d))$ of semistable (resp. stable) vector bundles of
rank $n$ and degree $d$ on $Y$. Let $J$ denote the generalised Jacobian
of degree $d$ on $Y$.

\begin{teo} {\bf (Theorem 3(A)).}\ \ {\it Let the
assumptions be above.

\begin{enumerate}
\renewcommand\labelenumi{{\rm (\alph{enumi})}}
\item {\rm Pic} $U^{\prime s} \approx $ {\rm Pic} $J \oplus \Z${\rm ,}

\item {\rm Pic} $U' \approx $ {\rm Pic} $J \oplus \Z${\rm ,}

\item $U'$ is locally factorial.
\end{enumerate}\vspace{-.4pc}

This completes the extension of results of {\rm \cite{6}} to nodal
curves.}\vspace{.5pc}
\end{teo}

Let $U = U(n,d)$ denote the moduli space of torsion-free sheaves of rank
$n$ and degree $d$ on $Y$. If $Y$ has only a single ordinary node as
singularity, then the variety $U(2,d)$ has a stratification, $U = U'
\cup U_1 \cup U_0$, a disjoint union. Points of $U_1$ correspond to
torsion-free sheaves $F$ of rank 2 with $F_y \approx {\cal O}_y \oplus
m_y$. Let $L$ be a rank $1$ torsion-free sheaf which is not locally
free. Let $U_{1,L}(2,d)$ be the subscheme of $U_1$ corresponding to
torsion-free sheaves of rank 2 with determinant isomorphic to $L$.

\begin{teo} {\bf (Theorem 2, Theorem 3(B)).}\ \ {\it Let $g_Y\ge2${\rm
;} if $g_Y = 2${\rm ,} assume that $d$ is odd for {\rm (b),} {\rm
(c),} {\rm (d)}. Then

\begin{enumerate}
\renewcommand\labelenumi{{\rm (\alph{enumi})}}
\item {\rm Pic} $U_{1,L}(2, d) \approx \Z ${\rm ,}\vspace{.1pc}

\item {\rm Pic} $U^s_1(2, d) \approx$ {\rm Pic} $J_X \oplus \Z${\rm ,}\vspace{.1pc}

\item {\rm Pic} $U_1(2, d) \approx$ {\rm Pic} $J_X \oplus \Z ${\rm ,}\vspace{.1pc}

\item $U_1(2,d)$ is locally factorial.
\end{enumerate}}
\end{teo}

In a subsequent paper, we study the Picard group of a seminormal
variety. As an application we compute the Picard groups of the
compactified Jacobian and some subvarieties of $U(2,d)$.

\begin{notation}
Let $Y$ denote an irreducible reduced projective curve with ordinary
nodes $y_j, j=1, \ldots, m$ as only singularities. Let $g$ be the
geometric genus and $g_Y$ the arithmetic genus of $Y$. For $y \in Y$,
let $({\cal O}_y,m_y)$ be the local ring at $y$. A torsion-free sheaf
$N$ on $Y$ is locally free on the subset $U$ of non-singular points of
$Y$. The rank $r(N)$ of $N$ is the\break rank of the locally free sheaf
$N\!\!\mid_U$. The degree $d(N)$ of $N$ is defined by $d(N) = \chi
(N)\break + r(N) (g-1)$, where $\chi$ denotes the Euler characteristic.
Let $N^*$ denote the torsion-free sheaf ${\rm Hom}(N,{\cal O})$.

Let $J$ and ${\overline J}$ be respectively the generalised Jacobian and
the compactified Jacobian of $Y$ (of a fixed degree) and ${\cal P} $ the
Poincar\'{e} bundle. Let $p_J$ denote the projection to ${\overline J}$.
Let $U=U(n,d)$ be the moduli space of semistable torsion-free sheaves of
rank $n$ and degree $d$ on $Y$. Let $\delta = gcd (n,d)$. Let $U'\subset
U$ be the open subvariety corresponding to vector bundles (i.e.
$S$-equivalence classes of $E$ such that ${\rm gr} E$ is a vector
bundle). Fix a rank $1$ torsion-free sheaf $L$ of degree $d$ on $Y$. Let
$U'_L $ (resp. $U_{1,L}$) be the subscheme of $U$ corresponding to
vector bundles (resp. torsion-free sheaves) with determinant isomorphic
to $L$ and $U_L$ its closure in $U$. Let $U^{\prime s} \subset U',
U^{\prime s}_{L} \subset U'_L $ etc. be the open subvarieties
corresponding to stable torsion-free sheaves. The variety $U$ is
seminormal (\cite{13}, Theorem~4.2), $U'$ and $U'_L$ are normal being
GIT-quotients of non-singular varieties \cite{10}. For $m=1$, $U$ has a
filtration $U \supset W_{n-1} \supset \cdots \supset W_0$, with $W_i$
seminormal closed subvarieties \cite{13}. $W_{i-1}$ is the non-normal
locus of $W_i, i= 1, \ldots ,n$ and $W_0$ is normal. Let $U'= U-W_1,
U_i= W_i-W_{i-1} (i=1,\ldots ,n-1), U_0=W_0$.
\end{notation}

\section{Torsion-free sheaves of rank 2}

In this section we study $U_L (2,d)$ and $U(2,d)$. Throughout the
section $E$ will denote a torsion-free sheaf of rank 2 and degree $d$ on
$Y$.

\begin{lem}
Let $E$ be a torsion-free sheaf with $\wedge^2 E=L$ torsion-free. Let
$N_1$ be a rank $1$ subsheaf of $E$ such that the quotient $N_2=E/N_1$ is
torsion-free.

\begin{enumerate}
\renewcommand\labelenumi{{\rm (\arabic{enumi})}}
\item If $N_1$ or $L$ is locally free{\rm ,} then $N_2 \approx
N_1^*\otimes L${\rm ,}\vspace{.1pc}

\item If $N_2$ is locally free{\rm ,} then $N_1 \otimes N_2 \approx L$.
\end{enumerate}\vspace{-.8pc}
\end{lem}

\begin{proof}
The canonical alternating form $E \times E \to L$ induces an
alternating form $N_1 \times N_1 \to L$. We claim that this form is
zero. This is clear at $y \in Y$ such that the stalk $(N_1)_y$ is free.
If $(N_1)_y \not\approx {\cal O}_y$, then $(N_1)_y = m_y$, also $L_y =
{\cal O}_y$ or $m_y$ (\cite{12}, Prop.~2, p.~164). Let $u,v$ be the two
generators of $(N_1)_y$. Since any ${\cal O}_y$-linear map from $m_y$
to $m_y$ (or ${\cal O}_y$) is given by the multiplication by $a \in
\overline{\cal O}_y$ (= normalisation of ${\cal O}_y)$ (\cite{12},
p.~169), the map $(N_1)_y \to L_y$ defined by $w \mapsto w \wedge u$ is
given by $w \wedge u = w a, a \in \overline{\cal O}_y$. In particular,
$0 = u \wedge u = u a$. Since $\overline{\cal O}_y$ is a domain, this
implies $a=0$. Thus $v \wedge u=0$ and hence $(N_1)_y \wedge\break (N_1)_y=0$.

Define an ${\cal O}$-bilinear map $b\!\!:\! N_1 \times N_2 \to L$ by $b(n_1,
n_2) = n_1 \wedge n_3$, where $n_3$ is a lift of $n_2$ in $E$. This is
well-defined as any two lifts $n_3, n'_3$ differ by an element of $N_1$
and $N_1 \wedge N_1 =0$ as seen above. The bilinear map $b$ induces an
injective sheaf homomorphism $N_2 \to {\rm Hom} (N_1, L)$ which is an
isomorphism outside the singular set of $Y$. If $N_1$ or $L$ is locally
free, then $d({\rm Hom} (N_1, L)) = d(L) - d(N_1)$ (\cite{4}, Lemma~2.5(B))
and hence $d({\rm Hom} (N_1, L)) = d(N_2)$. It follows that $N_2 \approx
{\rm Hom} (N_1, L)$.

If $N_2$ is locally free, the bilinear map $b$ gives an injective
homomorphism of torsion-free sheaves $N_1 \otimes N_2 \to L$. Since
$d(N_1 \otimes N_2)= d (N_1)+d (N_2) = d (L)$, this is an isomorphism.
This proves the lemma.

We remark that if both $N_1, N_2$ are not locally free then $N_1 \otimes
N_2$ has a torsion and $b$ gives a homomorphism $N_1 \otimes
N_2/$torsion $ \to L$ which is not an isomorphism.
\end{proof}

\begin{lem}
Assume that $Y$ has only one node $y$. Let $\pi\!\!:\! X \to Y$ be the
normalisation map and $\pi^{-1} y = \{x, z\}$. Let $N_1,N_2$ be line
bundles of degree $-1$ on $X$.

\begin{enumerate}
\renewcommand\labelenumi{{\rm (\alph{enumi})}}
\leftskip .05pc
\item Given a line bundle $L$ on $Y$ with $\pi^*L = N_1\otimes N_2 (x +
z)${\rm ,} there exists a vector bundle $E$ of rank $2$ and determinant
$L$ on $Y$ such that $E$ is $S$-equivalent to $\pi_*N_1 \oplus
\pi_*N_2$.

\item There exists a torsion-free sheaf $E$ of rank $2$ on $Y$ such that
$(1)$ $E_y \approx {\cal O}_y \oplus m_y${\rm ,} $(2)$ determinant of
$E$ is isomorphic to $\pi_* (N_1 \otimes N_2(z))$ and $(3)$ $E$ is
$S$-equivalent to $\pi_*N_1 \oplus \pi_* N_2$.
\end{enumerate}\vspace{-.8pc}
\end{lem}

\begin{proof}$\left.\right.$

\begin{enumerate}
\renewcommand\labelenumi{\rm (\alph{enumi})}
\leftskip .05pc
\item We shall construct a generalised parabolic bundle $(E', F_1(E'))$
on $X$ which gives the required vector bundle $E$ on $Y$. Take $E' =
L_1\oplus L_2, L_1 = N_1(x+z), L_2= N_2$.\break Let $e_1,e_2$ be basis
elements of $(L_1)_x, (L_1)_z$ respectively. Let $f_1,f_2$ be basis
elements of $(L_2)_x, (L_2)_z$ respectively. Define $F_1(E') = (e_2-f_1,
ce_1+f_2), c$ being a non-zero scalar. Since the projections $p_1, p_2$
from $ F_1(E')$ to $E'_x, E'_z$ are both isomorphisms, $E$ is a vector
bundle \cite{1}. Choose $c$ such that $L$ corresponds to the generalised
parabolic line bundle $(\pi^* L, (c,1)), (c,1)\in {\bf P}^1$ \cite{1}.
One has det$(E', F_1(E')) = ($det$~E', (c,1)) = (\pi^* L, (c,1))$. Hence
det~$E = L$. Since $F_1(L_1) = 0, \pi_* L_1(-x-z)$ is a sub-bundle of
$E$. The quotient is $\pi_* L_2$ as the projection from $F_1(E')$ to
$(L_2)_x \oplus (L_2)_z$ is onto. Thus $E$ is $S$-equivalent to
$\pi_*(N_1\oplus N_2)$.

\item Take $E'$ as in the above proof, define $F_1(E') = (e_1+f_2,f_1)$.
Since $p_1$ is an isomorphism and $p_2$ has rank $1$, $E_y \approx {\cal
O}_y \oplus m_y$. Since $(e_1+f_2)\wedge f_1 = 0 e_1\wedge e_2 +
f_1\wedge f_2 + \cdots,$ one has det$(E', F_1(E')) = (L_1\otimes L_2,
(0,1))$. Hence det$(E) = \pi_*(L_1\otimes L_2(-x)) = \pi_*(N_1\otimes
N_2(z))$. The final assertion follows as in the above proof.
\end{enumerate}\vspace{-1.4pc}
\end{proof}

\begin{propo}$\left.\right.$\vspace{.5pc}

\noindent Let $g_Y=1$. Then one has the following{\rm :}

\begin{enumerate}
\renewcommand\labelenumi{\rm (\arabic{enumi})}
\leftskip .05pc
\item $U_L (2,1) = \{$a point$\}$ for $L \in \overline{J}${\rm ,}\vspace{.2pc}

\item[] $ U(2,1) \approx \overline{J} \approx Y, U'(2,1)\approx J \approx
Y-\{$node$\}$.\vspace{.5pc}

\item $U_{\cal O} (2,0) \approx \overline{J}/i \approx \P^1${\rm ,} where
$i\!\!:\! \overline{J} \to \overline{J}$ is defined by $N \mapsto
N^*${\rm ,}\vspace{.2pc}

\item[] $U_L (2,0)  \approx \P^1$ and  $ U'_L (2,0) \approx \A^1${\rm ,} for $L\in J$.
\end{enumerate}
\end{propo}

\begin{proof}$\left.\right.$

\begin{enumerate}
\renewcommand\labelenumi{\rm (\arabic{enumi})}
\item For $y \in Y$, let $I_y$ denote the ideal sheaf of $y$. The dual
$I_y^*$ is a rank 1 torsion-free sheaf of degree 1 \cite{5}. It is
well-known that $ y \mapsto I_y^*$ gives an isomorphism $Y \to
\overline{J}^1$, where $\overline{J}^1$ is the compactified Jacobian of
degree 1 torsion-free sheaves.

\parindent 1pc Let $E$ be a stable rank 2 torsion-free sheaf of degree 1 on
$Y$. Then $h^1 (E)=0$ as $E$ is stable and hence $h^0 (E) =1$. Any non-zero section
$s \in H^0 (E)$ must be everywhere non-vanishing, otherwise it will
generate a rank 1 torsion-free subsheaf of degree $\geq 1$
contradicting the stability of $E$. Hence $s \in H^0 (E)$ generates a
unique trivial line sub-bundle ${\cal O}$ of $E$. The quotient $E / {\cal
O}$ must be torsion-free, if not then the kernel of $E \to (E / {\cal
O})/$torsion will contradict the stability of $E$. Thus we have a
morphism $h\!\!: U (2,1) \to \overline{J}^1$ given by $E \mapsto E/{\cal
O}$. Conversely, given $L \in \overline{J}^1$, Ext$^1 (L, {\cal O})=
H^1$ $(L^*)$ (\cite{4}, Proof of Lemma 2.5(B)). Since $h^0 (L^*)=0, h^1
(L^*)=1$, any non-zero element in Ext$^1 (L, {\cal O})$ determines
a unique (up to isomorphism) torsion-free rank 2 sheaf $E$ of degree 1.
It is easy to check that $E$ is stable. This gives the inverse of $h$.
Note that $h$ is in fact the determinant map.

\item We first prove that $W_0$ consists of a single point. Any element
in $W_0$ has stalk at the node $y$ isomorphic to $m_y \oplus m_y$. By
\cite{12}, Proposition~10, p.~174, such an element is the direct image
of a vector bundle $E_0$ on the desingularisation $\P^1$. Since $\pi_*
E_0$ is semistable, so is $E_0$. Hence $E_0 = {\cal O} (-1) \oplus {\cal
O} (-1)$. By Lemma 2.2(a), for every line bundle $L$ there exists a
vector bundle $E$ with determinant $L$ such that $E$ is $S$-equivalent
to $\pi_*({\cal O} (-1)) \oplus \pi_* ({\cal O} (-1))$. Thus for any $L
\in J$, $U_L (2,0)$ contains the point $\pi_* E_0$. One has $U_L \cap
W_1 = W_0$ \cite{1}. Thus every element of $U_L(2,0)$ is $S$-equivalent to
a vector bundle with determinant $L$. It follows that $U_L \approx
U_{\cal O}$.

We now prove that $U_{\cal O} (2,0) \approx \overline{J}/i \approx
\P^1$. Note first that the involution $i$ keeps the unique element
$\pi_* {\cal O}(-1)$ of $\overline{J}-J$ invariant and under the
isomorphism $Y \approx \overline{J}$, the map $\overline{J} \to
\overline{J}/i$ is the double cover $Y \to \P^1$ ramified at the image
of the node. Let $E$ be a semistable vector bundle of rank 2 with
trivial determinant. Let $E_1 $ be the vector bundle of degree 2
obtained by tensoring $E$ with a line bundle of degree 1. Since $E_1 $
is semistable, with slope $>0, h^1 (E_1) =0, h^0 (E_1) = 2$. Since the
evaluation map $Y \times H^0 (E_1) \to E_1$ cannot be an isomorphism,
there is a section of $E_1$ vanishing at a point and hence generating a
(torsion-free) subsheaf $N_1$ of rank 1, degree $\geq 1$. Since $E_1$ is
semistable, one must have $d(N_1)=1$. Hence $E$ has a rank 1 subsheaf
$N$ of degree $0$. The quotient $E/N$ is torsion-free in view of the
semistability of $E$.\break By Lemma~2.1(1), $E/N \approx N^*$. Thus $E$ is
$S$-equivalent to $N \oplus N^*$. Using the Poincar\'{e} bundle and the
properties of moduli spaces, one sees that this proves the proposition.
\end{enumerate}\vspace{-1.4pc}
\end{proof}

\begin{lem}
For $g_Y \geq 2, d$ even and $L \in J${\rm ,} one has
\begin{equation*}
{\rm codim}_{U'_L} (U'_L - U^{\prime s}_L) = 2g_Y - 3.
\end{equation*}
\end{lem}

\begin{proof}
A rank 2 vector bundle $E$ which is semistable but not stable contains
a torsion-free subsheaf $N_1$ with a torsion-free quotient $N_2 \approx$
Hom $(N_1, L) = N_1^* \otimes L$, where $L$ is determinant of $E$
(Lemma 2.1(1)). Thus $E$ is $S$-equivalent to $N_1 \oplus (N_1^*
\otimes L)$, hence dim $U'_L - U^{\prime s}_L = $ dim $J = g_Y$ and
codim$_{U'_L} U'_L - U^{\prime s}_L = 2g_Y -3 \geq 3$ if $g_Y \geq 3$.\vspace{.3pc}
\end{proof}

\begin{lem}$\left.\right.$

\begin{enumerate}
\renewcommand\labelenumi{{\rm (\arabic{enumi})}}
\item {\rm Codim}$_{U_L} (U_L-U'_L) \geq 3$ for $g_Y \geq 3$.\vspace{.1pc}

\item For $g_Y = 2${\rm ,} {\rm codim}$_{U_L} (U_L-U'_L) = 3$ if $d$ is
odd{\rm ,} $U_L = U'_L = {\P}^3$ if $d$ is even.
\end{enumerate}
\end{lem}

\begin{proof}$\left.\right.$

\begin{enumerate}
\renewcommand\labelenumi{\rm (\arabic{enumi})}
\item The points of $U_L - U'_L$ correspond to torsion-free sheaves
which are direct images of semistable vector bundles with fixed
determinant on partial normalisations of $Y$. Hence $U_L - U'_L$ is a
finite union of irreducible components each of dimension $3 (g_Y-1) -3$
$=3g_Y -6$ for $g_Y \geq 3$. Thus codim$_{U'_L} (U_L-U'_L) \geq 3$.

\item For $g_Y = 2$ the partial normalisations are of arithmetic genus
$1$. It follows from Proposition~2.3(1) that for $d$ odd, $U_L - U'_L$
consists of one or two points according as $g= 1$ or $g=0$. For $d$
even, $U_L = U'_L \approx {\P}^3$ (\cite{2}, Lemmas 3.3, 3.4, Corollary
3.5). We remark that Proposition 2.3(2) implies that the subset
$U_{0,L}$ of non-locally free sheaves in $U_L$ is isomorphic to ${\P}^1$
if $g=1$ and it consists of two smooth rational curves intersecting in a
point if $g=0$. The intersection point is the direct image of the unique
semistable bundle of degree $d-2$ on the desingularisation ${\P}^1$.
Note also that $U_{0,L} = U_L - U^{s}_L$ in this case.
\end{enumerate}\vspace{-1.5pc}
\end{proof}

\begin{lem}
{\rm Codim}$_{U'} U' - U^{\prime s} \geq 3$ for $g_Y \geq 3$ {\rm (}$d$
even{\rm )}.
\end{lem}

\begin{proof}
The surjective determinant map $U' \to J$ is a fibration with fibres
isomorphic to $U'_L, L$ a fixed line bundle of degree $d$. Hence the
lemma follows from Lemma 2.4.\vspace{.3pc}
\end{proof}

\begin{rem}$\left.\right.${\rm

\begin{enumerate}
\renewcommand\labelenumi{\rm (\arabic{enumi})}
\item Let $g_Y =1$. Then Pic $U(2,1) \approx G_m \oplus \Z$. For $L\in
J,$ Pic $U_L(2,0) \approx \Z$, and Pic $U'_L(2,0),$ Pic $U'(2,0),$ Pic
$U'(2,1)$ are trivial.\vspace{.15pc}

\item If $g_Y = 2$, then Pic $U'_L (2,d) \approx \Z \approx $ Pic
$U^{\prime s}_L(2,d)$ for all $d$.
\end{enumerate}}
\end{rem}

\begin{proof}
Part (1) follows from Proposition 2.3. Part (2) is proved in \cite{3}, \S2.4.
\end{proof}

\begin{propo}$\left.\right.$\vspace{.5pc}

\noindent For $g_Y \geq 3${\rm ,} one has{\rm :}

\begin{enumerate}
\renewcommand\labelenumi{\rm (\arabic{enumi})}
\item $U^{\prime s}_L (2,d) \approx \Z${\rm ,}\vspace{.15pc}

\item $U'_L (2,d) \approx \Z$.
\end{enumerate}\vspace{-.7pc}
\end{propo}

\begin{proof}
Let $p\!\!:\! \tilde{U}_L \to U_L$ be a (finite) normalisation. Since $U'_L$
is normal, $p$ is an isomorphism over $U'_L$ and $p$ gives a finite map
$\tilde {U}_L - p^{-1}U'_L \rightarrow U_L - U'_L$. Therefore codim
$\tilde{U}_L - p^{-1}U'_L=$ codim $U_L-U'_L \geq 3$ by Lemma~2.5. Since
$\tilde{U}_L$ is normal, this implies that Pic $\tilde{U}_L
\hookrightarrow \ {\rm Pic} (p^{-1} U'_L) \approx {\rm Pic} \ U'_L$.
Since $U_L$ is projective, so is $\tilde{U}_L$ and hence rank(Pic
$\tilde{U}_L ) \geq 1$. It follows that rank(Pic $U'_L) \geq 1$. Since
$U'_L$ is normal and by Lemma~2.4, codim$(U'_L -U^{\prime s}_L) \geq 3$
we have Pic $U'_L \hookrightarrow $ Pic $U^{\prime s}_L$. Thus rank(Pic
$U^{\prime s}_L) \geq 1$. By \cite{3}, Proposition~2.3, one has Pic
$U^{\prime s}_L \approx \Z $ or $\Z /m \Z, m \in \Z$. It follows that
Pic $U^{\prime s}_L \approx \Z$ and hence Pic $U'_L \approx \Z$.
\end{proof}

\begin{rem}{\rm
Putting together the results of \cite{3} and Proposition 2.8, we have
Theorem $1$.}
\end{rem}

\setcounter{subsection}{9}
\subsection{\it Varieties $U_1$ and $U_{1,L}$}

Henceforth we assume that there is only one node $y$. We first remark
that if $E$ is a rank 2 vector bundle then $E$ cannot be
$S$-equivalent to a direct sum of a line bundle and a non-locally free
torsion-free rank 1 sheaf. For, then, one has an exact sequence $0 \to
L_1 \to E \ \to L_2 \to 0$ with one of the $(L_1)_y $ or $(L_2)_y$
isomorphic to ${\cal O}_y$ and the other isomorphic to $m_y$. Since
${\rm Ext}^1 (m_y, {\cal O}_y)=0= {\rm Ext}^1 ({\cal O}_y, m_y)$, this
means $E_y \approx {\cal O}_y \oplus m_y$, i.e., $E$ is not locally free.
Similarly one sees that if $E_y \approx {\cal O}_y \oplus m_y$, then $E$
cannot be $S$-equivalent to a direct sum of two locally free sheaves. In
particular $E$ with $E_y \approx {\cal O}_y \oplus {\cal O}_y$ cannot be
$S$-equivalent to $E'$ with $E_{y}'$ not free unless
$[E]=[E'] \in W_0$. Hence taking determinant gives a well-defined
morphism det: $U'\cup U_1 \to \overline{J}_Y$ with det$(U')= J_Y$,
det$(U_1) = \overline{J}_Y-J_Y \approx J_{X}$. This morphism induces a
morphism of normalisations det: $P' \cup P_1 \to \tilde{J}_Y,
\tilde{J}_Y$ being the desingularisation of $\overline{J}_Y$ and $P',
P_1$ are respectively the pull backs of $U', U_1$ in the normalisation.

\setcounter{theore}{10}
\begin{lem}
Let $L \in \overline{J}_Y -J_Y$ with degree of $L$ even.

\begin{enumerate}
\renewcommand\labelenumi{{\rm (\arabic{enumi})}}
\item {\rm dim}$(U_{1,L} - U^s_{1,L}) = g_Y${\rm ,} for all $L${\rm
,}\vspace{.15pc}

\item {\rm codim} $U_{1,L} - U^s_{1,L} \geq 3$ for $g \geq 3$.
\end{enumerate}\vspace{-.7pc}
\end{lem}

\begin{proof}$\left.\right.$
\begin{enumerate}
\renewcommand\labelenumi{{\rm (\arabic{enumi})}}
\item From \S2.10, one sees that $E \in U_{1,L} -
U_{1,L}^s$ is $S$-equivalent to $N_1 \oplus N_2$ with one of $N_1,
N_2$ locally free and the other torsion-free but not locally free.
Also, one of them is a subsheaf and the other is a quotient sheaf.  By
Lemma 2.1, $E \sim M \oplus (M^* \otimes L), M \in J_Y$.  It follows
that dim$(U_{1,L} - U_{1,L}^s) =g_Y$.  In fact, one has $U_{1,L} -
U_{1,L}^s \approx J_Y$.

\item One has dim $U_{1,L}= 3g_Y-3$. Hence codim$(U_{1,L} - U_{1,L}^s)
= 2g_Y -3 \geq 3$ for $g_Y\break \geq 3$.
\end{enumerate}\vspace{-2pc}
\end{proof}

\pagebreak

\begin{lem}
For $L \in \overline{J}_Y-J_Y$ and $g_Y \geq 2${\rm ,} one has {\rm
codim}$_{U_L} (U_L - U_{1,L}) \geq 2$.
\end{lem}

\begin{proof}
The subset $U_L-U_{1,L}$ consists of torsion-free (semistable) rank 2
sheaves $E \approx \pi_* E_0, E_0$ semistable vector bundle of rank 2
on $X$ with det$E_0 \approx (\pi^* L/{\rm torsion} ) (- x )$ or $(\pi^*
L / {\rm torsion}) (-z)$ \cite{1}. Hence dim$(U_L - U_{1,L}) = 3g_X -3$
if $g_X \geq 2$, dim$(U_L - U_{1, L})=0$ if $g_X=1$ and $d$ is odd,
dim$(U_L- U_{1,L}) =1$ if $g_X=1$ and $d$ is even. Therefore, one has
for $g_Y \geq 3$, dim $U_L-U_{1,L} = 3g_Y-6$ and
codim$_{U_L}(U_L-U_{1,L}) = (3g_Y-3) - (3g_Y-6) =3$. For $g_Y=2$,
{\rm codim}$_{U_L}(U_L-U_{1,L})= 3$ if $d$ is odd and {\rm codim}$_{U_L}
(U_L-U_{1,L})=2$ if $d$ is even.
\end{proof}

\begin{lem}$\left.\right.$

\begin{enumerate}
\renewcommand\labelenumi{{\rm (\arabic{enumi})}}
\item $U_1^s$ is non-singular{\rm ,} $U_1$ is normal.\vspace{.15pc}

\item $U_{1,L}$ is normal{\rm ,} $U_{1,L}^s$ is non-singular.\vspace{.15pc}

\item $W_0^s$ is non-singular{\rm ,} $W_0$ is normal.
\end{enumerate}
\end{lem}

\begin{proof}
The moduli space $U$ is the geometric invariant theoretic quotient of
$R^{ss}$ by a projective linear group. Let ${\cal E} $ be the universal
quotient sheaf on $R^{ss} \times Y$. Let $R_1 = \{t \in R^{ss}\!\!\!\mid\!\!({\cal
E}_t)_y \approx {\cal O}_y \oplus m_y \}, R_0 = \{t \in
R^{ss}\!\!\!\mid\!\!({\cal E}_t)_y \approx m_y \oplus m_y \}$, $R_{1,L} = \{t \in
R_1\!\!\!\mid\!\!{\rm det} {\cal E}_t =L\}$. At any point $p \in R^{ss}$, the
analytic local model for $R_1 \hookrightarrow R^{ss}$ at $p$ is Spec
$A/(u,v) \hookrightarrow$ Spec $A$ where $A= \C[u,v] / (uv)$ (\cite{9},
Theorem~2(2), p.~576). Since the spectrum of a point is a regular
scheme, $R_1$ is regular. Since $U_1^s$ is a geometric quotient of
$R_1^s$, it follows that $U_1^s$ is a regular scheme. Since $R_1,
\overline{J}_Y - J_Y$ are regular and $R_{1,L}$ are all isomorphic,
$R_{1,L}$ is regular. Hence the assertion (2) follows. We remark here
that $R_1, R_{1,L}$ are not saturated for $S$-equivalence; $U_1$ and
$U_{1,L}$ are G.I.T. quotients of open subsets of $R_1$ and $R_{1,L}$
consisting of sheaves not $S$-equivalent to elements in $R_0$ and hence
are normal. The assertion (3) follows as (2) using \cite{9}, Theorem~2(3).\vspace{.3pc}
\end{proof}

\begin{propo}$\left.\right.$\vspace{.5pc}

\noindent Let $Y$ be an irreducible projective curve {\rm (}with one ordinary
node{\rm ),} $g_Y \geq 2$ and $n =2$.\break Then
\begin{equation*}
{\rm Pic} \ U_{1,L}^{s} \approx \Z \quad {\rm or} \quad \Z /m \Z, m
\in \Z.
\end{equation*}
\end{propo}

\begin{proof}
The idea of the proof is the same as that of \cite{6} or \cite{3},
Proposition~2.3. Hence we only indicate the necessary modifications. We
may assume $d \gg 0$. Then $R^1 p_{J *} {\cal P}^*$ is a vector bundle on
$\overline{J}_Y$. Let $\P = \P (R^1 p_{J_*} (\P^*)), \P_L=$ fibre of
$\P$ over $L \in \overline{J}_Y$. One has a universal family ${\cal E}$
of rank 2 torsion-free sheaves $E$ of degree $d$ on $\P \times Y$. Let
$\P^s, \P^s_L$ be the subvarieties corresponding to stable sheaves.
Since Ext$^1 ({\cal O}_y, {\cal O}_y) =0= $ Ext$^1 (m_y, {\cal O}_y)$,
one has $E_y \approx {\cal O}_y \oplus {\cal O}_y$ or ${\cal O}_y \oplus
m_y$. Hence by the universal property of moduli spaces, one has
morphisms $f_{\epsilon}\!\!: \P^s \to (U-W_o)^s$ and $f_{\epsilon, L}\!\!:
\P^s_L \to U^{\prime s}_L $ (or $U_{1,L}^s)$ if $L \in J_Y$ (or $L \in
\overline{J}_Y -J_Y)$. By \cite{10}, Chapter~7, Lemma~5.2$^{\prime}$, any
semistable torsion-free sheaf $E$ of $ d \gg 0$ is generated by global
sections. If $E_y \approx {\cal O}_y \oplus {\cal O}_y$ or ${\cal O}_y
\oplus m_y$, then by \cite{1}, Lemma~2.7, one has an exact sequence $0
\to {\cal O}_Y \to E \to G \to 0$ with $G$ torsion-free. Also $G
\approx$ det\,$E$ by Lemma~2.1(1). Hence $f_{\epsilon}$ and
$f_{\epsilon, L}$ are surjective. One shows that the induced map
$f^*_{\epsilon, L}$ on Picard groups is injective. This was checked in
\cite{3} for $L \in J_Y$, the same proof goes through for $L \in
\overline{J}_y -J_Y$ as $R_{1,L}^s$ and $U_{1,L}^s$ are \hbox{non-singular}
(Lemma~2.13(2)). Let $\P_{\overline{J}-J} = \P (R^1p_{J_*} ({\cal
P}^*\!\!\mid_{\overline{J}-J}))$ and $f_1\!\!: \P^s_{\overline{J}-J} \to U_1^s$. The
same argument gives that $f^*_1$ is injective and one has exact
sequences
\begin{align*}
&0 \to {\rm Pic}\, U_1^s \to {\rm Pic}\, \P^s_{\overline{J}-J} \to
\Z/{((n-1) d/a)\Z} \to 0, a = gcd (n,d),\\[.2pc]
&0 \to {\rm Pic}\, U^s_{1,L} \to {\rm Pic}\, \P^s_L \to \Z/{((n-1) d
/a) \Z} \to 0.
\end{align*}
Since $\P^s_L$ is an open subset of a projective space, Pic $\P^s_L$ is
isomorphic to $\Z$ or $\Z/{m \Z}$ and the same as true for Pic $U^s_{1,L}$.

We remark that the injectivity of $f^*_{\epsilon}$ does not seem to
follow similarly. In the notations of \cite{6}, Corollary 7.4, one
certainly gets a codimension one subvariety $\Gamma_0-\Gamma'_0$ of
$\Gamma_0$. Since $(U- W_0)^s$ is not necessarily non-singular it is not
clear that $\Gamma_0 - \Gamma'_0$ is a Cartier divisor, i.e., its ideal
sheaf is locally free. $U -W_0$ is seminormal, but not normal in
general, in particular it is not locally factorial.
\end{proof}

\begin{propo}$\left.\right.$\vspace{.5pc}

\noindent Let the notations be as in Proposition $2.14$. Then for $g_Y
\geq 3${\rm ,} $n=2$ and $g_Y=2${\rm ,} $n=2${\rm ,} $d$ odd{\rm ,} one has
\begin{equation*}
{\rm Pic}\, U_{1,L} \approx {\rm Pic}\, U_{1,L}^s
\approx \Z.
\end{equation*}
\end{propo}

\begin{proof}
For $d$ odd, $U_{1,L}= U^s_{1,L}$. Since $U_{1,L}$ is normal and
codim$(U_{1,L} - U^s_{1,L}) \geq 3$ (Lemma 2.11), Pic $U_{1,L}
\hookrightarrow$ Pic $U_{1,L}^s$ for $d$ even, $g_Y \geq 3$ as in the
proof of Proposition~2.8. Going to a finite normalisation we see that
rank (Pic $U_{1,L}) \geq 1$. We need Lemma~2.12 for this. The result now
follows from Proposition~2.14.
\end{proof}

\setcounter{subsection}{15}
\subsection{}

Assume that $g_Y =2, g_X=1, n=2, d=0$. Let $M$ be the moduli space of
$\alpha$-semistable GPBs $(E, F_1 (E))$ of rank 2, degree $0$ on a
smooth elliptic curve $X, 0 < \alpha <1, \alpha$ being close to $1$
\cite{1}. Let $M_L$ be the closed subscheme of $M$ corresponding to $E$
with determinant $L, L \in J_X$. Let $p_1\hbox{:}\ F_1 (E) \to E_x, p_2\hbox{:}\ F_1 (E)
\to E_z$ be the projections. Define $D_L = \{(E, F_1 (E)) \in
M_L|p_2$ has rank $\leq 1\}$ and $D_{1,L} = \{ (E,F_1 (E))
\in D_L|$rank $p_2=1, p_1$ isomorphism$\}$. $D_{1,L}$ is an open subscheme of
$D_L$ and $D_L$ is a closed subscheme of codimension 1 in $D$. There
is a surjective birational morphism $f\hbox{:}\ M \to U$ such that $D_L$ maps
onto $U_{L'}$ inducing an isomorphism $D_{1,L} \approx U_{1,L'}$ where $
L' = \pi_* (L(-z))$. We shall determine $D_L, D_{1,L}$ explicitly and
use the explicit description to compute Pic~$U_{1,L'}$. Note that $D_L
\approx D_{\cal O}$ for all $L$.

\setcounter{theore}{16}
\begin{propo}$\left.\right.$\vspace{.5pc}

\noindent $D_L$ is isomorphic to a $\P^2$-bundle over $\P^1$. Outside $\P^1 - \{4$
points$\}${\rm ,} this bundle is of the form $\P ({\cal O} \oplus
\epsilon)${\rm ,} $\epsilon$ being a rank $2$ vector bundle.
\end{propo}

\begin{proof}
It is not difficult to check that $(E,F_1 (E))$ of degree $0$, rank 2 is
$\alpha$-semistable if and only if $E$ is a semistable vector bundle and
for any line sub-bundle $L$ of $E$ of degree $0, F_1 (E) \neq L_x \oplus
L_z$. Moreover, $(E, F_1 (E))$ is $\alpha$-stable if and only if $E$ is
semistable and $F_1 (E) \cap (L_x \oplus L_z) =0$ for any sub-bundle of
degree 0.
\pagebreak

Let $e_1,e_2$ and $e_3, e_4$ be the bases of $E_x$ and $E_z$
respectively. The subspace $F_1 (E)$ defines a point in the Grassmannian
${\rm Gr}$ of two-dimensional subspaces of $V=E_x \oplus E_z$. Let ${\rm Gr} \subset
\P (\wedge^2V)$ be the Pl\"ucker embedding, let $(X_1, Y_1, X_2, Y_2,
X_3,Y_3)$ be the Pl\"{u}cker coordinates. Any element in $\wedge^2 V$ is of
the form $X_1 e_1 \wedge e_2+ Y_1 e_3 \wedge e_4+ X_2 e_1 \wedge e_4+
Y_2 e_2 \wedge e_3+ X_3 e_3 \wedge e_1+Y_3 e_2 \wedge e_4$. The
Grassmannian quadric is given by $X_1 Y_1+X_2 Y_2+X_3 Y_3=0$. Since $E$
is semistable, one has either (a) $E=M \oplus M^*, M \in J_X$ or (b)
there is a non-trivial extension $0 \to M_1 \stackrel{g}{\to} E
\stackrel{h}{\to} M_2 \to 0$ with $M_1 \approx M_2 \approx M \in J_X,
M^2 = {\cal O}$. In either case $E$ is an extension of $M_2$ by $M_1;
M_1, M_2 \in J_X$. Choose $e_1, e_2, e_3, e_4$ to be basis elements of
$(M_1)_x ,(M_2)_x, (M_1)_z, (M_2)_z$ respectively. Let $ D_V\subset {\rm Gr}$
be defined by\break $Y_1 = 0$.

\begin{case}{\rm
Assume that $E = M_1 \oplus M_2, M_1^* = M_2, M_1 \neq M_2$. The group
$\P$ (Aut $E)= \P (G_m \times G_m) \approx G_m$ acts on $D_V \subset \P
(\wedge^2 V)$ by $t (X_1, X_2, Y_2, X_3,Y_3) = (X_1, X_2, Y_2, t X_3,
t^{-1} Y_3)$. It is easy to see that $D_V /\!/ G_m \approx \P^2$, the
quotient map $D_V \to \P^2$ being given by $(X_1, X_2, Y_2, X_3, Y_3)
\mapsto (X_1,X_2, Y_2)$. Let $D_{1,V} = D_V - \{ (X_1=0)$ $\cup
(1,0,0,0,0)\}$. The image of $D_{1,V}$ in $\P^2$ is given by $\P^2 -
\{{(X_1=0)} \cup (1,0,0)\}$.

Let ${\cal P}_X \to J_X \times X$ be the Poincar\'{e} bundle, ${\cal P}_x
= {\cal P}|_{J'_X \times x}, {\cal P}_z = {\cal P}|_{J'_X
\times z}, J'_X = J_X - J_2$, $J_2$ being the group of 2-torsion points of
$J_X$. The group $G_m \times G_m$ acts on the bundles $\V = ({\cal P}_x
\oplus {\cal P}_x^*) \oplus ({\cal P}_z \oplus {\cal P}^*_z)$, and
$\wedge^2 \V$ as above, giving $G_m$-action on $\P (\wedge^2 \V)$ and
$D_{\V} /\!/ G_m \approx \P^2$-bundle over $J'_X$. This $\P^2$-bundle is
in fact the bundle $\P ({\cal O} \oplus ({\cal P}_x \otimes {\cal
P}^*_z) \oplus ({\cal P}_z \otimes {\cal P}_x^*))$. The involution on
$J_X$ given by $i (M) = M^*$ lifts to an action on this bundle
(switching second and third factors), hence it descends to a bundle on
$J'_{X} / i = \P^1 - \{ 4$ points$\}$, of the form $\P ({\cal O} \oplus
\epsilon), \epsilon$ a vector bundle of rank 2 on $J'_X/i$.}
\end{case}

\begin{case}{\rm
There are, up to isomorphism, exactly four bundles $E$ given by extension
of type (b). Since any automorphism of $E$ is of the form $\lambda Id +
\mu g \circ h$, one has $\P $ (Aut $E) \approx G_a$ under the
isomorphism $(\lambda, \mu) \mapsto t = \mu \lambda^{-1} \in G_a$. The
action of $G_a$ on $V$ is given by $te_1 =e_1, te_3= e_3, te_2= e_2+
te_1, te_4 = e_4 + te_3$ and that on $D_V$ is given by $t(X_1,X_2, Y_2, X_3,
Y_3) = (X_1, Y_2+ tY_3, X_2+ tY_3, X_3- t(X_2+Y_2) - t^2 Y_3, Y_3)$. It
is not difficult to see that the ring of invariants for $G_a$-action on
$D_V$ (resp. on the hyperplane $Y_1=0 $ of $\P (\wedge^2 V))$ is
generated by $X_1, X_2-Y_2, Y_3$ (resp. $X_1, X_2-Y_2, Y_3, X_2Y_2+X_3
Y_3)$. The non-semistable points for the $G_a$-action are
$\{X_1=Y_3=X_2-Y_2 =0\}$. It follows that $D_V /\!/ G_a \approx \P^2$, the
quotient map $D_V \to \P^2$ being given by $(X_1, X_2, Y_2, X_3, Y_3)
\to (X_1, X_2-Y_2, Y_3)$. Clearly, $D_{1,V} /\!/ G_a \approx \P^2 - (\{X_1
=0\} \cup (1,0,0))$. We remark that non-stable GPBs correspond to the
line $Y_3=0$ in $\P^2$. In case $E=M_1 \oplus M_2,$\break $M_1 = M_2$ with
$M_1^2={\cal O}$, one sees that corresponding quotient $D_V /\!/ G_a$ is
$\P^1$ which is identified to the line $Y_3=0$ in the above $\P^2$. Note
that there are no stable GPBs in the last case.

It follows that there is a $\P^2$-fibration $\phi\hbox{:}\ D_L \to \P^1$ which
is locally trivial outside the set of four points in $\P^1$. By Tsen's
theorem (\cite{8}, p.~108, Case (d))$, \phi$ is a locally trivial
fibration. This completes the proof.}\vspace{-.5pc}
\end{case}
\end{proof}

\begin{coro}$\left.\right.$\vspace{.5pc}

\noindent Let $g_X=1${\rm ,} $g_Y=2${\rm ,} $d$ even{\rm ,} $n=2$.

\begin{enumerate}
\renewcommand\labelenumi{\rm (\arabic{enumi})}
\item $U_{1,L'}$ is non-singular.\vspace{.15pc}

\item {\rm Pic} $U_{1,L'} \approx \Z$.
\end{enumerate}
\end{coro}

\begin{proof}$\left.\right.$

\begin{enumerate}
\renewcommand\labelenumi{\rm (\arabic{enumi})}
\item It follows immediately from the proof of Proposition 2.18 that
$D_{1,L}$ is a (locally trivial) fibration over $\P^1$ with non-singular
fibres isomorphic to $\P^2- \{(X_1=0) \cup (1,0,0)\}$. Hence $ D_{1,L}$
and $U_{1,L'}$ are non-singular.

\item $D_L - D_{1,L} \cong $ (hyperplane $H) \cup \{$a line $\ell \}, H
\cap \ell = \Phi$, Pic $D_L \approx $ Pic $\P^1 \oplus $ Pic $\P^2$.
Since $D_L$ is non-singular, $0 \to \Z H \to $ Pic $D_L \to $
Pic$(D_L-H) \to 0$ is exact. It follows that Pic $D_L -H \cong $ Pic $\P^1
= \Z$. Since $\ell$ is of codimension 2, Pic$(D_{1,L}) \approx $
Pic$(D_L-H) \cong \Z$. Thus Pic $U_{1,L'}\approx $ Pic $D_{1,L} \approx
\Z$.
\end{enumerate}\vspace{-1.7pc}
\end{proof}

\begin{rem}{\rm
Note that $H \to \P^1$ is a $\P^1$-bundle. The fibres of this bundle are
given by $X_1=0$ in $D_V$, the restriction of this bundle to $\P^1 - \{
4$ points$\}$ is $\P (\epsilon)$. Under the map $D_L \to U_{L'}$, this
$\P^1$-bundle maps onto one component in $U_{L'} - U_{1,L'} $
isomorphic to $ J_X / i~(\approx \P^1)$. This component corresponds to
sheaves of the form $\pi_* E_0$, det $E_0 \approx L(-x-z)$. The line
$\ell$ maps isomorphically onto the other component isomorphic to
$\P^1$, it corresponds to $\pi_* E_0$, det $E_0 \approx L(-2z)$. Since
$g_X =1$, $E_0$ are semistable but not stable. Thus unlike in the case
when $L'$ is a line bundle ($Y$ smooth or nodal) $U_{L'} - U^s_{L'}$ is
not the Kummer variety. It has an open subset isomorphic to $J_Y$ (Proof
of Lemma~2.11(1)) whose complement is the union of two disjoint smooth
rational curves.

Putting together Proposition 2.15 and Corollary 2.18, we have proved the
following.}
\end{rem}

\begin{therr}[\!]
Let $Y$ be an irreducible projective curve of arithmetic genus $\geq 2$
with only a single ordinary node as singularity. Let $L$ be a rank $1$
torsion-free sheaf which is not locally free. Then
\begin{equation*}
{\rm Pic}\ U_{1,L} \approx \Z.
\end{equation*}
\end{therr}\vspace{-1pc}

\section{Pic and local factoriality of $\U\,'\,\hbox{\bf (}\n\hbox{\bf
,}\,\d\hbox{\bf ),}\ \U_{{\bf 1}\hbox{\co}\,\La}\,\hbox{\bf (2}\hbox{\bf
,}\,\d\hbox{\bf )}$}

\subsection{}

In this section we prove Theorems 3A and 3B. Throughout the section, we
assume that $n \ge 2$ and if $n \ge 3$ then $g\ge 2$. One has a map
$U'_L \times J \rightarrow U'$ given by tensorisation. We first remark
that Pic $U'$ cannot be computed easily using this map. The map induces
a map of Picard groups Pic $U' \approx $ Pic $U'_L\, \oplus $ Pic $J
\rightarrow $ Pic $U'_L\, \oplus $ Pic $J$. The induced map Pic $J
\rightarrow $ Pic $J$ is not identity, it is multiplication by $n$. The
right map to consider is the determinant morphism which does induce
identity on Pic $J$ as we show below:

\begin{teo} {\bf 3A.}\ \ {\it One has the following{\rm :}

\begin{enumerate}
\renewcommand\labelenumi{{\rm (\alph{enumi})}}
\item {\rm Pic} $U^{\prime s} \approx$ {\rm Pic} $J \oplus \Z${\rm
,}\vspace{.1pc}

\item {\rm Pic} $U' \approx$ {\rm Pic} $J \oplus \Z ${\rm ,}\vspace{.1pc}

\item $U'$ is locally factorial.
\end{enumerate}}\vspace{-1pc}
\end{teo}

\begin{proof}$\left.\right.$

\begin{enumerate}
\renewcommand\labelenumi{\rm (\alph{enumi})}
\item Without loss of generality, we may assume that $d \gg 0$. Then a
semistable vector bundle $E$ of degree $d$ is globally generated
(\cite{10}, Lemma~5.2) and contains a trivial sub-bundle of rank $n-1$.
Let $\P =\P (R^{1}_{p_{j_*}} ({\cal P}^* \otimes \C^{n-1}))$, it is a
projective bundle over $J$. Let $\P_L$ denote its fibre over $L \in J,
\P_L$ is a projective space. $\P$ parametrises a family ${\cal E}$ of
vector bundles on $Y$ of rank $n$, degree $d$ and containing a trivial
sub-bundle of rank $n-1$. Let $\P^s =\{p \in \P|{\cal E}_p$
stable$\}, \P^{s}_{L} = \P^s \cap \P_L$. One has canonical surjective
morphisms $f\hbox{:}\ \P^s \rightarrow U^{\prime s} (n,d), f_L\hbox{:}\ \P^{s}_{L}
\rightarrow U^{\prime s}_{L}(n,d)$ such that the induced maps
$f^*$: Pic $U^{\prime s} \rightarrow $ Pic $\P^s, f^{*}_{L}\hbox{:}$ Pic $U^{\prime s}_L
\rightarrow$ Pic $\P^{s}_{L}$ are injective (\cite{3}, Proposition~2.3;
\cite{6}, Propositions~7.6, 7.8, 7.9). Clearly, Pic $\P \approx $ Pic
$J\, \times $ Pic $\P_L \approx $ Pic $J \times \Z$. Under the conditions of
the theorem we know that (\cite{3}, Theorem I) Pic $U^{\prime s}_{L}
\approx \Z $ and hence Pic $\P_{L}^{s} \approx \Z$. Hence the surjective
restriction map Pic $\P_L \rightarrow $ Pic $\P_{L}^{s}$ is an
isomorphism for all $L \in J$. Hence codim$_{\P_L}(\P_L - \P_{L}^{s})
\neq 1$ and therefore codim$_{\P} ( \P-\P^s) \geq 2$. Thus Pic $\P^s
\approx$ Pic $\P \approx$ Pic $J \oplus \Z$ and hence
\begin{equation*}
\hskip -.55cm {\rm Pic} \ U^{\prime s} \hookrightarrow {\rm Pic} \ J
\oplus \Z.
\end{equation*}

The natural map $p\hbox{:}\ \P^s \rightarrow J$ factors as $p=$ det
$\circ f$, where det is the determinant map $E \mapsto
\stackrel{n}{\wedge} E$. Since both $f$ and det are surjections, so is
$p$. Note that $f^* \circ {\rm det}^* =p^*\hbox{:}$ Pic $J \rightarrow $
Pic $\P^s$ is injective. It follows that ${\rm det}^*$ is injective.

\parindent 1pc One has the following diagram with the last column exact.
\begin{align*}
\hskip -.55cm \begin{array}{lllllllll}
0 && 0 && 0 \\ [2.5mm]
\downarrow && \downarrow && \downarrow \\ [2.5mm]
{\rm Pic}\,J &=&{\rm Pic}\,J &=& {\rm Pic}\,J \\ [2.5mm]
\downarrow && \downarrow && \downarrow \\ [2.5mm]
{\rm Pic}\,U' &\rightarrow & {\rm Pic}\,U^{\prime s}& \hookrightarrow& {\rm Pic}\,J \oplus \Z \\
[2.5mm]
\downarrow && \downarrow && \downarrow \\ [2.5mm]
{\rm Pic}\,U'_L & \stackrel{\approx\ \,}{\rightarrow}& {\rm Pic}\,U^{\prime s}_{L} &
\stackrel{\approx\ \,}{\rightarrow} & \Z \\ [2.5mm]
\downarrow && \downarrow && \downarrow \\ [2.5mm]
0&& 0& &0
\end{array}
\end{align*}
Here $\Z$ denotes the image of Pic $U^{\prime s}_{L}$ in Pic
$\P^{s}_{L}$. The map Pic $U^{\prime s} \rightarrow $ Pic $U^{\prime
s}_{L}$ is the restriction map and is surjective (\cite{3},
Proposition~3.2 and 3.5). It now follows from the diagram that the
injection Pic $U^{\prime s} \rightarrow $ Pic $J \oplus \Z$ is an
isomorphism and the second column is exact.

\item and (c). Since codim$_{U'}(U' - U^{\prime s}) \geq 2$ under the
conditions of the theorem and $U'$ is normal (\cite{3}, Proposition~3.4(i)),
it follows that the restriction map Pic $U' \rightarrow $ Pic
$U^{\prime s}$ is injective. The restriction morphism Pic $U'
\rightarrow $ Pic $U'_L$ is surjective (\cite{3}, Propositions~3.2, 3.5). The
restriction map Pic $U'_L \rightarrow $ Pic $U^{\prime s}_{L}$ is an
isomorphism \cite{3}. It now follows from the commutative diagram that
Pic $ U' \approx $ Pic $U^{\prime s}$ under the restriction map. By
arguments similar to those in the proof of \cite{3}, Proposition~3.6,
this implies that $U'$ is locally factorial.
\end{enumerate}\vspace{-1.6pc}
\end{proof}

\begin{teo} {\bf 3B.}\ \ {\it
Let $Y$ be an irreducible projective curve of arithmetic genus $g_Y\geq
2$ with only a single ordinary node as singularity. If $g_Y = 2${\rm ,}
then assume that $d$ is odd. Let $L$ be a rank $1$ torsion-free sheaf
of degree $d$ which is not locally free. Let $U_{1,L}$ be the subscheme
of $U$ corresponding to torsion-free sheaves of rank $2$ with
determinant isomorphic to $L$.

\begin{enumerate}
\renewcommand\labelenumi{{\rm (\alph{enumi})}}
\item {\rm Pic} $U^s_1 \approx$ {\rm Pic} $J_X \oplus \Z${\rm ,}\vspace{.2pc}

\item {\rm Pic} $U_1 \approx$ {\rm Pic} $J_X \oplus \Z ${\rm ,}\vspace{.2pc}

\item $U_1$ is locally factorial.
\end{enumerate}}
\end{teo}

\begin{proof}
The proof is more or less identical with that of Theorem~3A. One has
only to\break replace $f, f_L$ by the maps $f_1, f_{\epsilon,L}$ of
Proposition~2.14 and use Theorem~2 instead of Theorem~1.
\end{proof}

\section{The dualising sheaves of $\U\,'$ and $\U\,'_{\La}$}

\subsection{}

Let $K(Y)$ denote the Grothendi\'{e}ck group of vector bundles on $Y$. Then
$K(Y) \approx \Z\, \oplus$\, Pic $Y$ under the map $[E] \mapsto$ (rank $E$,
det $E), [E]$ being the class of a vector bundle $E$ in $K(Y)$. The
inverse map is given by $ n \mapsto [n \cdot {\cal O}_Y]$ for $n \in \Z$
and $L \mapsto [L] - [{\cal O}_X]$ for $L \in $ Pic $Y$.

Let $\chi = d+n (1-g), P(m)= \chi +rm$, fix $m \gg 0$. Let $Q= {\rm
Quot}(\C^{P(m)} \otimes {\cal O}_Y (-m), P)$ be the Hilbert scheme (`the
Quot scheme') of quotients of $\C^{P(m)} \otimes {\cal O}_Y (-m)$ with
Hilbert polynomial $P$. Let ${\cal F} \rightarrow Q \times Y$ be the
universal family. Let $R_m \subset Q$ be the open subset consisting of
$q \in Q$ such that $H^1 ({\cal F}_q (m)) =0, H^0 (\sum (m)) \simeq H^0
({\cal F}(m))$ under the canonical map, $\sum = \C^{P(m)} \otimes {\cal
O}_Y (-m)$. The open subvariety $R^{ss}$ of $Q$ consisting of $q \in Q$
such that ${\cal F}_q$ is a semistable torsion-free sheaf is contained
in $R_m$. The subset $R^{\prime ss}$ of $R^{ss}$ corresponding to semistable
vector bundles is a smooth variety, so is the closed subset $R^{\prime ss}_{L}
\subset R^{\prime ss}$ consisting of semistable vector bundles with fixed
determinant $L$ (\cite{10}, Remark, p.~167).

The moduli space $U'$ (resp. $U'_L)$ is a geometric invariant theoretic
good quotient of the smooth irreducible scheme $R^{\prime ss}$ (resp.
$R^{\prime ss}_{L})$ by the group $G = P$(Aut $\sum) \approx PGL (N), N
\gg 0$ \cite{10,12}. The restriction of the universal family on $Q \times
Y$ gives a universal family ${\cal F} \rightarrow R^{\prime ss}_{L}
\times Y$ of vector bundles on $Y$ of rank $n$, degree $d$. Let Pic\,$^G
(R^{\prime ss}_{L})$ denote the group of line bundles on $R^{\prime
ss}_{L}$ with $G$-action (compatible with the $G$-action on $R^{\prime
ss}_{L})$. For a vector bundle $E$ on $Y$, one defines an element
$\lambda_{\cal F} (E) \in {\rm Pic}^G (R^{\prime ss}_{L})$ by
\begin{equation*}
\lambda_{{\cal F}} (E):= \otimes_i ({\rm det} R^{i}_{p_{1*}} ({\cal F}
\otimes p^*_2 E))^{(-1)^{i+1}},
\end{equation*}
where $p_1$ and $p_2$ are projections to $R^{\prime ss}_{L}$ and $Y$
respectively. $ \lambda_{{\cal F}} (E)$ depends only on the class of $E$
and $\lambda_{{\cal F}}\hbox{:}\ K (Y) \rightarrow {\rm Pic}^G (R_{L}^{\prime ss})$ is a
group homomorphism.

\setcounter{theore}{1}
\begin{propo}$\left.\right.$\vspace{.5pc}

\noindent Let $E$ be a vector bundle on $Y$ with {\rm rank}$(E) = n /
\delta${\rm ,} {\rm det} $(E) = {\cal O}_Y (- \frac{\chi}{\delta})${\rm
,} $\chi = d+n$ $(1-g)${\rm ,} $\delta = gcd (n,d)$. Then
$\lambda_{{\cal F}} (E)$ descends to $U'_L (n,d)$ as the generator $\L$
of {\rm Pic} $U'_L (n,d)$.
\end{propo}

\begin{proof}
By \cite{3}, Propositions~3.2, 3.5, the generator $\L$ is obtained by
the descent of the line bundle $\L'$ on $R^{\prime ss}_{L}$ given by
\begin{equation*}
\L'= ({\rm det}\ Rp_{1_*} {\cal F})^{\frac{n}{\delta}} \otimes
(\stackrel{n}{\wedge} ({\cal F}|_{R^{\prime ss}_{L} \times y_0}))^{\chi /
\delta},
\end{equation*}
$y_0$ being a non-singular point of $Y$. Here det $Rp_{1_*} {\cal F}$
denotes the determinant of cohomology (\cite{7}, Ch.VI, pp.~135--136).
However, our definition is different from the standard one, it is the
inverse of the line bundle defined in \cite{7} as det $Rp_{1_*} {\cal
F}$. One has det $R_{p_{1*}} ({\cal F}) = \lambda_{\cal F}$ $(1), 1=$
class of ${\cal O}_{Y}$. If $h$ denotes the class of the structure sheaf
of the point $y_0, h=[{\cal O}_Y (y_0)]$ $- [{\cal O}_Y]$, then we claim
that
\begin{equation*}
\stackrel{n}{\wedge} {\cal F} \mid_{R_{L}^{\prime ss} \times y_0} =  -
\lambda_{\cal F} (h).
\end{equation*}
\end{proof}\vspace{-1pc}

\begin{potc}
For $m \gg 0$ one has the exact sequence
\begin{equation*}
0 \rightarrow {\cal F} (m) \rightarrow {\cal F} (m+1) \rightarrow {\cal
F} (m) \mid_{R^{\prime ss}_{L} \times y_0} \rightarrow 0,
\end{equation*}
${\cal F} (m) = {\cal F} \otimes {\cal O}_Y (m), {\cal O}_{Y} (1)$ being
a line bundle of degree 1 on $Y$. Since $R^{1}_{p_{1_*}} ({\cal F} (m'))
=0$ $\forall m' \geq m, R^1p_{1_*} ({\cal F} (m)\!\! \mid_{R^{\prime
ss}_{L} \times y_0}) =0$, the direct image sequence gives
\begin{equation*}
0 \rightarrow R^0p_{1_*} ({\cal F} (m)) \rightarrow R^0p_{1_*}
({\cal F} (m+1)) \rightarrow R^0p_{1_*} ({\cal F}
(m) |_{R^{'ss}_{R_L} \times y_0}) \rightarrow 0.
\end{equation*}

$\left.\right.$\vspace{-1.5pc}

\noindent Since det $p_{1_*} ({\cal F} (m')) = - \lambda _{{\cal F}}
(1 + m'h), m' \geq m$, and
\begin{equation*}
{\rm det} (p_{1_*} {\cal F}(m) |_{R^{\prime ss}_{L} \times y_0}) \approx
{\rm det} (p_{1_*} {\cal F}|_{R^{\prime ss}_{L} \times y_0} ) =
\stackrel{n}{\wedge} {\cal F} \mid_{R^{'ss}_{L} \times y_0},
\end{equation*}
one has
\begin{align*}
\stackrel{n}{\wedge} {\cal F} |_{R^{\prime ss}_{L} \times y_0}
&=- \lambda_{\cal F} ((m+1)h) + \lambda_{\cal F} (1+mh)\\
&=- \lambda_{\cal F} (h).
\end{align*}
This proves the claim.

Thus we have
\begin{align*}
\L' &= \frac{n}{\delta} \lambda_{\cal F} (1) - \frac{\chi}{\delta}
\lambda_{\cal F} (h)\\[2mm]
&= \lambda_{{\cal F}} \left(\frac{n}{\delta} - \frac{\chi
h}{\delta}\right) = \lambda_{\cal F} (E).
\end{align*}
\end{potc}\vspace{-1pc}

\begin{rem}{\rm
Note that the line bundle $\L'$ exists on $R^{ss}$ and descends to $U'$
(\cite{3}, Proposition~3.5). Also $\lambda_{\cal F} (E)$ makes sense for
${\cal F} \rightarrow R^{ss} \times Y$, the universal family on $R^{ss}
\times Y$. The above relation between $\lambda_{\cal F} (E)$ and $\L \in $
{\rm Pic} $U'_L (n,d)$ holds for $\lambda_{\cal F} (E)$ and $\L \in $
{\rm Pic} $U'(n,d) \approx {\rm Pic}\, U'_L \oplus {\rm Pic}\, J$.}
\end{rem}\vspace{.1pc}

\setcounter{subsection}{3}
\subsection{\it Computation of the dualising sheaves}

Both $U'$ and $U'_{L}$ are normal and Cohen--Macaulay as they are
quotients of smooth varieties by $PGL(N)$. They are also locally
factorial (\cite{3}, Theorem~2; Theorem~1). A locally factorial
Cohen--Macaulay variety is Gorenstein, i.e., its dualising sheaf
$\omega$ is locally free. The tangent sheaf $T_{U'}$ of $U'$ is locally
free on the smooth open subscheme $U^{\prime s}$ of codimension $\geq
2$. Hence the determinant of $T_{U'}$ defines a line bundle det $T_{U'}$
on $U'$. Since it coincides with $\omega^{-1}$ on $U^{\prime s}$, it
follows that $\omega^{-1} = {\rm det}\ T_{U'}$. Similarly one has a
locally free dualising sheaf $\omega_{L}$ on $U'_L$ with
$\omega^{-1}_{L} = {\rm det}\ T_{U'_L}$.

\setcounter{theorr}{3}
\begin{therr}[\!]
Let the assumptions be as in Theorem~{\rm 1}. Then one has the following{\rm :}
\begin{enumerate}
\renewcommand\labelenumi{{\rm (\alph{enumi})}}
\item $\omega {\approx} -2 \delta \L${\rm ,} $\L=$ generator of {\rm Pic} $U'_L
(n,d)${\rm ,}

\item Let $F_0$ be a vector bundle on $Y$ of rank $2r$ and degree
$2(-d+r(g-1))$.  Then $\omega \approx \lambda_{\cal F} (F_{\circ})
\otimes$ {\rm det} $\wedge${\rm ,} where $\wedge$ is a line bundle on $J$ given
by
\begin{equation*}
\hskip -.5cm \wedge = {\rm det} (p_{J!} [{\cal P}] \otimes {\rm det}
p_{J!} [{\cal P}^*])^{r-1} \otimes {\rm det} p_{J!} ([{\cal P} \otimes
p^{*}_{2} F_0])^{-1}.
\end{equation*}
\end{enumerate}
\end{therr}

\begin{proof}
In view of the injective morphism $f^{*}_{L}\hbox{:}\ {\rm Pic}\,U'_L
\rightarrow {\rm Pic}\, P^{s}_{L}$ mapping $\L$ to ${\cal O}_{p^{s}_{L}}
(\frac{d}{\delta} (r-1))$, it suffices to prove that
\begin{equation*}
{\rm det} f^*_L T_{U'_L} \approx {\cal O}_{\P^{s}_{L}} (2d (r-1)).
\end{equation*}
One has $f^* T_{U'} \approx R^{1}_{p_{\P^{s}_{*}}} ({\cal E}^* \otimes
{\cal E}), f^{*}_{L} T_{U'_{L}} \approx R^{1}_{p_{\P^{s}_{L*}}}\!\!\!(Ad
{\cal E} ) \approx R^{1}_{p_{\P^{s}_{*}}}\!\!\,(Ad {\cal
E})\!\!\mid_{\P^{s}_{L}}$. Also, det $R^{1}_{p_{\P^{s}_{*}}} ({\cal E}^*
\otimes {\cal E}) \approx {\rm det}\ R^{1}_{p_{\P^{s}_{*}}} (Ad {\cal
E})$, so that det $f^{*}_{L} T_{U'_{L}} \approx $ det $
R^{1}_{p_{\P^{s}_{*}}} ({\cal E} \otimes {\cal
E}^*)\!\!\mid\!\!\P^{s}_{L}$.
\end{proof}

\subsection*{\it Computation of {\rm det} $R^{1}_{p_{\P^{s}_{*}}} ({\cal
E} \otimes {\cal E}^*)$}

There is a universal exact sequence on $\P^s \times Y$.
\begin{equation}
0 \rightarrow {\cal O}_{\P^s \times Y} \otimes \C^{r-1} \rightarrow
{\cal E} \rightarrow (1 \times p)^* {\cal P} \otimes p^*_{\P^s} {\cal
O}_{\P^s} (-1) \rightarrow 0.
\end{equation}
For $d \gg 0, H^0 ({\cal E}^{*}_{t}) =0 \forall t \in \P^s, H^0 ({\cal
E}_t \otimes {\cal E}^{*}_{t})$ consists of scalars as ${\cal E}_t$ is
stable. Hence by tensoring (1) with ${\cal E}^*$ and taking direct
images, one gets (for $d \gg 0$ and $ (1 \times p)^* =p^{\#})$
\begin{align*}
0 \rightarrow {\cal O}_{\P^s} &\rightarrow {\cal O}_{\P^{s}} (-1)
\otimes p_{\P^{s}_{*}} (p^{\#} {\cal P} \otimes {\cal E}^*) \rightarrow
R^{1}_{p_{\P^{s}_{*}}} ({\cal E}^* \otimes \C^{r-1})\\
&\rightarrow R^{1}_{p_{\P^{s}_{*}}} ({\cal E}^* \otimes {\cal E})
\rightarrow  0.
\end{align*}
Hence,
\begin{align}
{\rm det}\, R^{1}_{p_{\P_{*}^{s}}} ({\cal E} \otimes {\cal E}^*) &\approx
{\rm det} (R^{1}_{p_{\P^{s}_{*}}} {\cal E}^*)^{r-1}\nonumber\\
& \otimes {\rm det} ({\cal O}_{\P^s} (-1) \otimes p_{\P^{s}_{*}} p^{\#}
{\cal P} \otimes {\cal E}^*)^{-1}.
\end{align}
$R^{1}_{p_{\P^{s}_{*}}} ({\cal E}^*)$ is computed by taking dual of (1)
and direct images as follows:\vspace{.6pc}

\hskip 3.2pc $0 \rightarrow p^{\#} {\cal P}^* \ \otimes \ p^{*}_{\P^{s}}
{\cal O}_{\P^s} (1) \rightarrow {\cal E}^{*} \rightarrow {\cal O}_{\P^s
\times Y} \ \otimes \ \C^{r-1} \rightarrow 0.$\hskip 1.8cm $(1)^{*}$\\

\noindent Since $p_{\P^{s}_{*}} p^{\#} {\cal P}^*=0 =p_{\P^s*} ({\cal
E}^*)$ for $d \gg 0$, one has the direct image sequence
\begin{align*}
0 &\rightarrow {\cal O}_{\P^s} \otimes \C^{r-1} \rightarrow {\cal
O}_{\P^s} (1) \otimes R^{1}_{p_{\P^{s}_{*}}} (p^{\#} {\cal P}^*) \rightarrow
R^{1}_{p_{\P^{s}_{*}}} {\cal E}^*\\
&\rightarrow {\cal O}_{\P^s} \otimes \C^{(r-1)g} \rightarrow 0
\end{align*}
and hence
\begin{equation*}
{\rm det}\, R^{1}_{p_{\P^{s}_{*}}} {\cal E}^* \approx {\rm det} ({\cal O}_{\P^s} (1)
\otimes R^{1}_{p_{\P^{s}_{*}}}(p^{\#} {\cal P}^*)).
\end{equation*}
Since $h^1 ({\cal P}^{*}_{t}) = - \chi ({\cal P}^{*}_{t}) = d+g-1$ for
$t \in J$, one gets
\begin{equation}
{\rm det}\, R^{1}_{p_{\P^{s}_{*}}} {\cal E}^* \approx {\cal O}_{\P^s} (d+g-1)
\otimes {\rm det}\, R^{1}_{p_{\P^{s}_{*}}} (p^{\#} {\cal P}^*).
\end{equation}
Tensoring $(1)^*$ with $p^{\#} {\cal P}$ gives
\begin{equation*}
0 \rightarrow p^{*}_{\P^s}{\cal O}_{\P^s} (1) \rightarrow {\cal
E}^* \otimes p^{\#} {\cal P} \rightarrow {\cal O}_{\P^s \times Y}
\otimes \C^{r-1} \otimes p^{\#}
{\cal P} \rightarrow 0,
\end{equation*}
and hence the direct image sequence
\begin{equation*}
0 \rightarrow {\cal O}_{\P^s} (1) \rightarrow p_{\P^{s}_{*}} ({\cal E}^*
\otimes p^{\#} {\cal P}) \rightarrow p_{\P^{s}_{*}}(\C^{r-1} \otimes
p^{\#} {\cal P}) \rightarrow 0.
\end{equation*}

By tensoring with ${\cal O}_{\P^{s}} (-1)$ and taking det, one has
\begin{equation*}
{\rm det} (p_{\P^{s}_{*}} (p^{\#} {\cal P} \otimes {\cal E}^*) \otimes {\cal
O}_{\P^s} (-1)) \approx {\rm det} (p_{\P^{s}_{*}} (p^{\#} {\cal P} \otimes
\C^{r-1}) \otimes {\cal O}_{\P^s} (-1)).
\end{equation*}

$\left.\right.$\vspace{-1.5pc}

\noindent Since $h^0 ({\cal P}_t) = d+1-g $ for $t \in J$, the latter is
isomorphic to det $p_{\P^{s}_{*}} (p^{\#} {\cal P} \otimes \C^{r-1})
\otimes {\cal O}_{\P^{s}} ((g-d-1) (r-1))$. Thus we have
\begin{align}
&{\rm det} (p_{\P^{s}_{*}} (p^{\#} {\cal P} \otimes {\cal E}^{*}) \otimes
{\cal O}_{\P^{s}} (-1))\nonumber\\
&\quad\ \approx {\rm det}\, p_{\P^{s}_{*}} (p^{\#} {\cal P} \otimes
\C^{r-1}) \otimes {\cal O}_{\P^{s}} ((r-1) (g-d-1)).
\end{align}
Substituting in (2) from (3) and (4) gives
\begin{equation}
{\rm det}\, R^{1}_{p_{\P^{s}_{*}}} ({\cal E}^* \otimes {\cal E}) \approx
{\cal O}_{\P^{s}} (2 (r-1)d) \otimes \Delta^{r-1},
\end{equation}
where $\Delta^{-1} = {\rm det} (R^{1}_{p_{\P^{s}_{*}}} p^{\#} {\cal P}^*)
\otimes {\rm det} (p_{\P^{s}_{*}} p^{\#} {\cal P})$.

Since $\Delta |_{\P_{L}^{s}}$ is trivial, from (5) one has
\begin{equation*}
{\rm det} f^{*}_{L} T_{U'_L} \approx {\cal O}_{\P^{s}_{L}} (2(r-1)d),
\end{equation*}
this proves (a).

If $F_{0}$ is a vector bundle of rank $2r$ and degree $2(-d+r(g-1))$,
then from sequence (1), one sees that
\begin{equation*}
\lambda_{{\cal E}} ([F_0]) \approx {\cal O}_{\P^s} (-2d (r-1)) \otimes
{\rm det}^* (p_{J_*} {\cal P} \otimes p_{Y}^{*} F_0),
\end{equation*}
so that (5) becomes
\begin{equation*}
{\rm det} (R^1p_{\P^{s}_{*}} ({\cal E} \otimes {\cal E}^{*})) \approx
\lambda_{{\cal E}} ([F_0])^{-1} \otimes {\rm det}^* (p_{J_*}({\cal P}
\otimes p^{*}_{Y} F_0)) \otimes \Delta^{r-1}.
\end{equation*}
Since $p = {\rm det} \circ f, p^* = f^* \circ {\rm det}^*$ and $f^*$ is
injective, (b) also follows.

\section*{Acknowledgement}

We would like to thank C~S~Seshadri for useful discussions
and correspondence.

\end{document}